\documentclass[12pt,leqno]{article}
\usepackage{amsmath,amsthm,amsfonts,amssymb,amscd}
\usepackage[latin1]{inputenc}
\usepackage{graphics}
\usepackage{graphicx}

\title{\bf Extension of germs of holomorphic foliations}
\date{\today}
\author{Gabriel Calsamiglia and Paulo Sad}
\begin{document}

\maketitle

\begin{abstract}
We consider the problem of extending germs of plane holomorphic foliations to foliations of compact surfaces. We show that the germs that become regular after a single blow up and admit  meromorphic first integrals can be extended, after local changes of coordinates, to foliations of  compact surfaces. We also show that the simplest elements in this class can be defined by polynomial equations. On the other hand we prove that, in the absence of meromorphic first integrals there  are uncountably many elements without polynomial representations.
\end{abstract}

\section {\bf Introduction}

In this paper we treat the following problem: let ${\mathcal Fol}(\mathbb C^2,0)$ be the set of germs of holomorphic foliations defined in a neighborhood of $0\in \mathbb C^2$ which are singular at the origin; we consider two such foliations to be equivalent when they are conjugated by a local holomorphic diffeomorphism of $\mathbb C^2$ at $0\in \mathbb C^2$. We select some family $L\subset {\mathcal Fol}(\mathbb C^2,0)$ and ask if an equivalence class contains a foliation that is defined by polynomial differential equations. Any member of such a class admits an extension to a foliation of the complex projective plane after a suitable local change of coordinates.

A very simple example is given by the set of hyperbolic singularities in the Poincaré domain, namely, foliations defined by the 1-form
$$
(x+A(x,y))\,dy - (\lambda y +B(x,y))\,dx=0,
$$
where $\lambda \notin \mathbb R$ and $A$ and $B$ are holomorphic functions such that $A(0,0)=B(0,0)=0$ and whose derivatives at $0\in \mathbb C$ vanish. By the theorem of linearisation of Poincaré (\cite {PO}), we have holomorphic equivalence to the linear part
$$
x\,dy - \lambda y\,dx = 0.
$$
It is interesting to notice that when $\lambda \le 0$ but $\lambda \notin \mathbb Q$ it is not known if  any equivalence class contains a foliation defined by polynomial equation (see \cite {YO}).

In their paper \cite {G-T}, Genzmer and Teyssier introduced a tool that allows to treat this problem from an analytic point of view. Roughly speaking the idea is: given the family $L$, a surjective map $\psi$ from the set $[L]$ of equivalence clases in $L$ to a space $I$ of invariants is defined; it is assumed that there are appropriate topologies to turn $\psi$ into an "analytic" map. The image of the equivalence classes of polynomially defined foliations is then a countable union of analytic manifolds (of finite dimension!) and cannot be the whole of $I$ provided $I$ is a "huge" space. This is a sort of analytic Baire property. As an application, the authors consider the family $L$ of saddle-node singularities of Milnor number 2 (in fact the choice of the Milnor number is not relevant). Those are singularities defined by forms
$$
[y(1+\mu x)+ R(x,y)]\,dx - x^2\,dy=0
$$
where ${\rm ord}_{(0,0)}(R) \ge 3$. According to \cite {M-R}, these singularities can be obtained applying a  convenient gluing procedure to normal forms
$$
y(1+\mu x)\,dx - x^2\,dy=0\,;
$$
in this gluing process non trivial local holomorphic diffeomorphisms $h(z)=z+\dots$ are used and become elements of the space of invariants $I$. The conclusion is that there are uncountably many saddle-nodes which are not equivalent to saddle-nodes defined by polynomial equations.

In the present paper we apply these ideas to the space $T$ of singular foliations that become regular after a single blow-up. The lack of singularities implies that every leaf of the foliation intersects the exceptional divisor and a finite number of them are tangent to it. In order to simplify the exposition, we consider the case with just one point of tangency of order 1 with the exceptional divisor; let $T_1$ be the space of these singularities. We may then state

\vspace{0.2in}
\noindent {\bf Theorem 1.} There are foliations in $T_1$ which are not holomorphically equivalent to foliations defined by polynomial equations.

\vspace{0.2in}
\noindent The space of invariants in this situation contains the local involutions at the tangency points which are naturally associated to the foliations in $T_1$.

\vspace {0.3in}
On the other hand, there are some situations where  we may wonder if the foliation  is conjugated to a polynomial model. For example, let us consider a germ of holomorphic function $f:(\mathbb C^2,0)\rightarrow (\mathbb C,0)$ such that $0\in \mathbb C^2$ is an isolated singularity of $f$. By a theorem of Mather, $f$ has finite determinacy: there exists $k\in \mathbb N$ such that the $k$-jet $f_k$ of $f$ at $0\in \mathbb C^2$ is conjugated to $f$ by means of a holomorphic diffeomorphism $\phi$, i.e., $f= {f_k}\circ \phi$.  It follows that the foliation defined by $d\,f=0$ is conjugated by $\phi$ to the foliation defined by $d\,f_k=0$. The question of whether a similar statement holds true for foliations defined by meromorphic functions arises: is a foliation in $(\mathbb C^2,0)$ defined by meromorphic function equivalent to a polynomially defined foliation? A theorem by Cerveau and Mattei (\cite {C-M}) gives sufficient conditions on the function to conclude that it is the case: let $f/g$ be a germ at $0\in \mathbb C^2$ of a meromorphic function ($f$ and $g$ ate supposed to be relatively prime germs of holomorphic functions) such that the 1-form $fd\,g-gd\,f$ defines a foliation with an isolated singularity at $0\in \mathbb C^2$. Then $f/g$ has finite determinacy, that is, $f/g$ is conjugated to $f_k/g_k$ for some $k\in \mathbb N$, where $f_k$ and $g_k$ are the k-jets of $f$ and $g$ at $0\in \mathbb C^2$. The simplest examples of functions defining foliations in $T_1$ do not satisfy the hypothesis of Cerveau and Mattei's result - take for instance $\frac {y^2-x^3}{x^2}=const$. However we prove the following

\vspace {0.2in}
\noindent {\bf Theorem 2.} Any foliation in $T_1$ admitting a generic meromorphic first integral is equivalent to one given by polynomial equations.
\vspace{0.2in}

The genericity condition can be stated as follows: a meromorphic first integral exists that vanishes with multiplicity one on the unique separatrix that is tangent to the exceptional divisor. This theorem improves the result in \cite {CA} in the generic case. In \cite{CA},  Casale proves that any germ admitting a meromorphic first integral is equivalent to a germ of a foliation on an algebraic surface (not necessarily the projective plane). By using the results in \cite{CAL} we are able to extend his result to foliations in $T$ defined by meromorphic functions. In order to state the theorem along the same lines of \cite{CA}, we use the notion of  {\bf algebraic-like foliation}: we say that a germ of holomorphic foliation $\mathcal F$ is an {\bf algebraic-like foliation} when there exists an algebraic surface $S$ and an algebraic foliation of $S$ which is equivalent to $\mathcal F$ in a neighborhood of some singularity.

\vspace {0.2in}
\noindent {\bf Theorem 3.} Any foliation in $T$ admitting a meromorphic first integral is an algebraic-like foliation.

\vspace {0.2in}

%\noindent We remark two things: 1) the simplest possible example of a  foliation in $T$ with a meromorphic first integral, the one defined as $\frac %{y^2-x^3}{x^2}=const$, does not satisfy the hypothesis of the theorem of Cerveau and Mattei. 2) the definition of "algebraic-like" foliation comes from the %paper of Casaleï¿½s , where Theorem 2 is proven for the case of a foliation in $T_1$. Actually we give an improvement for this result showing

A good point to discuss is whether Theorem 3 can be stated replacing "algebraic-like" foliation by "polynomially defined" foliation. We do not know the answer in general.

Sections 2, 3 e 4 of the present paper are devoted to the proof of Theorem 1, and Sections 5 and 6 to prove  Theorems 2 and 3. The two blocks are independent and can be read separately.

We thank G. Smith for useful conversations and especially for the example of the quintic in section \ref{sec:T1}.

\section {\bf Preliminaries}

We consider the set $T_1$ of germs of holomorphic foliations in $(\mathbb C^2,0)$ defined by differential 1-forms of the type

\begin {equation} \label{def1}
\sum_{j\ge 2}b_j(x,y)dx - \sum_{j\ge 2}a_j(x,y)dy = 0
\end{equation}

\noindent where $a_2(x,y)=xy, b_2(x,y)=y^2$ and $xb_3(x,y)-ya_3(x,y)=\beta x^4, \beta \ne 0$.

\vspace{0.2in}
After one blow-up $(x,t)\mapsto(x,tx)$, the  foliation is regular, with only one point of tangency of order 1 with the exceptional divisor (the equation is normalized as to have the tangency point given by $t=0$).

\vspace {0.2in} To each ${\cal F}\in T_1$ we can associate a local involution $i_{\cal F}(t)$ defined for $t\in \mathbb C$ close to $0\in \mathbb C$; moreover, it can be easily seen that for a holomorphic family $\alpha\in U\subset {\mathbb C^m}\mapsto {\cal F_{\alpha}}\in T_1$, the function $(\alpha,t)\mapsto i_{{\cal F}_{\alpha}}(t)$ is holomorphic.

\vspace {0.2in}
Let $Inv:= \{i(t)= \sum_{j\ge 1}a_jt^j \in {\mathbb C}\{t\}, a_1=-1, i\circ i(t)=t\}$; we consider in $\mathbb C\{t\}$ the norm $||\sum_{j\ge 0}c_jt^j||:= \sum_{j\ge0}\dfrac{|c_j|}{j!}$, which induces a distance $d$. Since

\begin {equation} \label{def2}
Inv_k :=\{i(t)\in \mathbb C \{t\}; i(0)=0, i^{\prime}(0)=-1 \,\,{\rm and}\,\, i\circ i(t)=t \,\,{\rm mod}\,\, t^{k+1}\}
\end {equation}

\noindent is closed in $(\mathbb C\{t\},d)$ for each $k\ge 1$, and $Inv= \cap_{k\ge 1}Inv_k$, we conclude that $Inv$ is closed in $(\mathbb C\{t\},d)$.

\vspace {0.2in}
Now we take

\begin{equation} \label{def3}
{\cal L}_1 \{t\} := \{\sum_{j\ge 0}c_jt^j \in \mathbb C\{t\};\,\, \sum_{j\ge 0}|c_j|<\infty \}
\end{equation}

\noindent Clearly ${\cal L}_1 \{t\}$ is a vector subspace of $\mathbb C\{t\}$; any power series in ${\cal L}_1\{t\}$ defines a holomorphic function whose domain of definition    contains the   unit disc $ {\mathbb D}= \{z\in \mathbb C; \,\, |z|\le 1\}$. On the other hand, the Taylor series centered at $0\in \mathbb C$ of a holomorphic function  defined in a neighborhood of $\bar {\mathbb D}$ belongs to ${\cal L}_1 \{t\}$.

\vspace {0.2in}
We define $||\sum_{j\ge 0}c_jt^j||_1:= \sum_{j\ge 0}|c_j|$ for $\sum_{j\ge 0}c_jt^j \in {\cal L}_1\{t\}$; with this norm ${\cal L}_1\{t\}$ becomes a Banach space. Let $d_1$ be the associated distance.

\vspace {0.2in}
\noindent {\bf Lemma 1.} The inclusion map from $({\cal L}_1\{t\},d_1)$ to $({\mathbb C}\{t\},d)$ is continuous.
\proof It is enough to remark that

\begin {equation} \label{continui}
||\sum_{j\ge 0}c_jt^j|| = \sum_{j\ge 0} \dfrac {|c_j|}{j!} \leq \sum_{j\ge 0}|c_j|= ||\sum_{j\ge 0}c_jt^j||_1
\end {equation}
\qed

\vspace {0.2in}

It follows that $Inv\cap {\cal L}_1\{t\}$ is closed in $({\cal L}_1\{t\},d)$. Therefore, $Inv\cap {\cal L}_1\{t\}$, endowed with the metric $d_1$, becomes a complete metric space, in particular a Baire space.

\section{\bf Realizing Involutions}

We introduced in the last section a map $i$  that takes foliations of $T_1$ to involutions of $\mathbb C\{t\}$.

\vspace {0.2in}
\noindent {\bf Lemma 2.} The map $i: T_1 \longrightarrow Inv$ is surjective.
\proof 1) Given some $i(t)\in Inv$, we construct first a local foliation around the disc ${ \mathbb D}\times \{0\}$ which has a tangency point at $(0,0)$ with this disc and whose associated involution is $i(t)$. We start  by mapping ${\mathbb D}\times 0$ to ${\mathbb C}\times 0$ via some holomorphic diffeomorphism $\phi$ which satisfies
\begin {itemize}
\item $\phi (0)= 0$.
\item $\phi$ conjugates \,$t\mapsto -t$\, to $i(t)$.
\end {itemize}
We then extend $\phi$ to some holomorphic diffeomorphism $\Phi$ in a neighborhood of ${\mathbb D}\times \{0\}$ and define the foliation $\cal H$ as the image by $\Phi$ of the foliation defined as $d(x-t^2)=0$.

\noindent 2) The next step consists in the following gluing process:
\begin {itemize}
\item we take the surface $S$ obtained after blowing-up $\mathbb C^2$ at $(0,0)$, foliated by $dt=0$ ($(x,y)$ are coordinates in $\mathbb C$, $(t=\dfrac {y}{x},x)$ are coordinates in $S$). In $S$ we remove a disc $\mathbb D{_{\frac {1}{2}}} \times U$, where $U$ is a small neighborhood of $0\in \mathbb C$.
\item the trivial foliation $dt=0$ in $\overline {({\mathbb D}\setminus {\mathbb D}_{\frac{1}{2}})}\times U$ is equivalent to the restriction of $\cal H$ to a region $R$ given as a local saturation (along the leaves of $\cal H$) of some annulus ${\mathbb A}\times \{0\}$
around $(0,0)\in {\mathbb D}\times \{0\}$. This equivalence is then used to glue ${\cal G}|_{ S\setminus \{(\overline {{\mathbb D}\setminus {\mathbb D_{\frac{1}{2}}}}) \times U \}}$ with $\cal H$; since it can be taken close to the Identity, the resulting foliation is defined around a $(-1)$-curve and is thus equivalent to the blow-up of an element of $T_1$. \qed
\end{itemize}

\section{Adapting Genzmer-Teyssier}
\label{adapt}
Our aim is to show that there are foliations in $T_1$ which are not holomorphically equivalent to any foliation in $T_1$ defined by a polynomial equation. In order to do that, we
need to change  the map $i$. Let $G$ be the group of Moebius  transformations of $\mathbb P^1$ which fix $0\in \mathbb C$ (in the t-coordinate associated to the blow up). We consider the map

\begin{equation} \label{grupo}
I: G\times T_1 \longrightarrow Inv,\,\,\, I(g,{\cal F})= g^{-1}\circ i_{\cal F}\circ g .
\end {equation}

\vspace{0.2 in}
\noindent {\bf Remark:} In fact the map of Lemma 2 induces a bijection between $[T_1]$ and $Inv/G$ (see \cite{CAL}).
\vspace{0.2 in}

Let $T_1^{(k)}$ denote the subset of elements of $T_1$ defined by a polynomial equation of degree $k$. The goal is therefore to prove that

\begin {equation} \label {goal}
\cup_k I(G\times T_1^{(k)}) \ne Inv
\end {equation}

We follow the procedure exposed in \cite{G-T}. We have to prove that the image of an embedding $\xi: {\overline {\mathbb D}}^l\longrightarrow (Inv,d)$ leaves a trace in $Inv \cap {\cal L}_1\{t\}$ which has empty interior in the topology defined by $d_1$.

Let us consider then some $f\in Im(\xi)\cap {\cal L}_1\{t\}$ and $0<\lambda<1$. Any power series defined as $f_{\lambda}(t)= {\lambda}^{-1}f({\lambda}t)$ belongs to ${\cal L}_1\{t\}$ and $d_1(f_{\lambda},f)\rightarrow 0$ as $\lambda\rightarrow 1$; furthermore, the radius of convergence of $f_{\lambda}$ is greater than $1$. If for some sequence ${\lambda}_m \rightarrow 1$ it happens that  $f_{{\lambda}_m}\notin Im(\xi)$, we are done; otherwise we replace $f$ by some  $d_1$-close $f_{\bar \lambda}$  and we still have $f_{\bar \lambda}\in Im(\xi)\cap {\cal L}_1\{t\}$. In order to simplify the notation we  use $f$ instead of $f_{\bar \lambda}$.

We then have $f= -t+\sum c_jt^j \in Im(\xi)\cap {\cal L}_1\{t\}$, with radius of convergence greater than 1. The tangent space $T_{f}Im(\xi)$ has some finite dimension $l$. Any element in $T_{f}Im(\xi)$ is a power series $\sum a_jt^j \in \mathbb C\{t\}$; after truncating the elements of $T_{f}Im(\xi)$ up to some sufficiently high order $m_0$, we still have a linear subspace of dimension $l$. Therefore, for each $m\ge m_0$, a power series in $T_{f}Im(\xi)$ is completely determined once we know the first $m$ coefficients.

Now we consider the path $\alpha (u):=h_u^{-1}\circ f\circ h_u$, where $h_u(t)= t+ut^m$ for $m\ge 0$. Clearly $h_u^{-1}$ is well defined in some disc of radius greater than 1 for $|u|$ small enough. This guarantees that $\alpha (u)$ is inside $Inv \cap {\cal L}_1\{t\}$. The tangent vector $\alpha^{\prime}(0)$ ({\it which we intend to prove that is transverse to} $T_fIm(\xi)$) has its $(m-1)$-jet equal to zero, therefore $\alpha^{\prime}(0)=0$ if it belongs to $T_f Im(\xi)$. But an easy computation shows that
\begin{equation} \label {path}
\alpha (u)(t)= h_u^{-1}\circ f\circ h_u (t) = -t + \sum_{j=2}^{m-1}c_jt^j + (c_m-2u)t^m + \cdots
\end{equation}
and then
\begin{equation} \label {contra}
\alpha^{\prime}(0)= -2t^m + \cdots
\end{equation}
which is a contradiction that proves Theorem 1.

\section {A Model}

Let us consider a holomorphic foliation $\mathcal G$ defined in some open set of $\mathbb C^2$ which contains $0 \in \mathbb C^2$ as an isolated singularity; we assume that the exceptional divisor is not invariant for the blown-up foliation. We will now conjugate $\mathcal G$ to a special model $\mathcal F$.

\begin{itemize}

\item {\bf Step 1}: we blow up at $0\in \mathbb C^2$; the exceptional divisor $E_1$ is not invariant for the blown-up foliation $\tilde {{\mathcal G}_1}$. We select a point $p\in E_1$ where  $\tilde {{\mathcal G}_1}$ is transverse to $E_1$ and take a neighborhood $V_1$ of this point where  $\tilde {{\mathcal G}_1}$ is trivial. In parallel, we blow up at $0\in \mathbb C^2$ the trivial foliation $d\,y=0$ to a foliation $\tilde {\mathcal G}_2$ which now has the exceptional divisor $E_2$ as an invariant set (with one singularity). We take a regular point of $\tilde {\mathcal G}_2$ in $E_2$ and a neighborhood $V_2$ of this point where $\tilde {\mathcal G}_2$ is trivial. We then glue $\tilde {{\mathcal G}_1}$ to $\tilde {\mathcal G}_2$ by a holomorphic diffeomorphism from $V_1$ to $V_2$ which sends ${\tilde {{\mathcal G}_1}}|_{V_1}$ to ${\tilde {\mathcal G}_2}|_{V_2}$. We get a surface which contains two divisors, still denoted  by $E_1$ and $E_2$, with $E_1 \cdot E_1=E_2\cdot E_2= -1$ and $E_1\cdot E_2=1$, and a foliation $\tilde {\mathcal G}$ conjugated to $\tilde {{\mathcal G}_1}$ and $\tilde {\mathcal G}_2$ in neighborhoods of $E_1$ and $E_2$ respectively.

\item {\bf Step 2}: we consider now the surface obtained after blowing up ${\mathbb D}\times {\mathbb P}^1$ at some point of $\{0\}\times \mathbb P^1$; we have inside it two divisors  $E_1'$ and $E_2'$ such that $E_1'\cdot E_1'=E_2'\cdot E_2'=-1$ and $E_1' \cdot E_2'=1$. Since a neighborhood of $E_1 \cup E_2$ is biholomorphically equivalent to a neighborhood of $E_1' \cup E_2'$ by a diffeomorphism that takes $E_1$ to $E_1'$ and $E_2$ to $E_2'$, we may define a foliation $\tilde {\mathcal F}$ in a neighborhood of $E_1' \cup E_2'$ as the image of $\tilde {\mathcal G}$. The blow-down of the restriction of  $\tilde {\mathcal F}$  to a neighborhood of $E_1'$ is the model $\mathcal F$ we mentioned above. We may see  $\tilde {\mathcal F}$ as the blow-up at some point of $\{0\}\times \mathbb P^1$ of a foliation ${\mathcal F}_1$ defined in ${\mathbb D}\times {\mathbb P}^1$.
\end{itemize}

\noindent In other words, modulo holomorphic equivalence,  $\mathcal G$ is  obtained by blowing-up
a foliation defined in ${\mathbb D}\times {\mathbb P}^1$ at some point of transversality with
$\{0\}\times \mathbb P^1$. As can be checked there are many choices involved in the construction and the model obtained in the product is not unique at all. Our strategy to find extensions to algebraic foliations will be to extend these types of foliated products to some foliation on $S\times\mathbb{P}^1$ where $S$ is a compact Riemann surface.

\section {Algebraic Case}

Let us remark that if $\mathcal G \in T$     then $\mathcal F_1$ is regular along $\{0\}\times \mathbb P^1$ as well (we are keeping the same notation used in the last Section). Furthermore, if $\mathcal G$ has a meromorphic first integral then
the same is true for ${\mathcal F}_1$; in particular a first integral $R(x,t)$ can be seen as a holomorphic family of rational functions $x\in {\mathbb D} \mapsto R_x(t)=R(x,t)\in {\mathbb P}^1$ of some degree $d$. It is not difficult to show that $x\in {\mathbb D} \mapsto R_x$ is locally injective at $x=0$.

To prove Theorem 2 we are going to see that up to reparametrizing the $x$-variable, we can extend this family to one parametrized by $\mathbb{P}^1$. The induced foliation will then define an extension of $\mathcal{F}_1$ to $\mathbb{P}^1\times\mathbb{P}^1$. To prove Theorem 3 what we are going to do is approximate this family by one passing through $R_0(t)$ and whose parameter space is a compact Riemann surface $S$. The associated foliation on $S\times\mathbb{P}^1$ --an algebraic surface-- can be blown up at an adequate point of $\{0\}\times {\mathbb P}^1$ so that the obtained foliation approximates $\tilde {\mathcal F}$ and will satisfy the necessary conditions to be  conjugated  to $\tilde {\mathcal F}$ in a neighborhood of $E_1'$ (according to \cite {CAL}).

\subsection{Proof of Theorem 2}
\label{sec:T1}
The proof of Theorem 2 does not use approximation and can be done after a suitable change of the first integral and of the coordinates on the product $\mathbb{D}\times\mathbb{P}^1$. In this subsection we suppose that $\mathcal{G}$ is a foliation in $T_1$ admitting a meromorphic first integral and consider its model $\mathcal{F}_1$ and a first integral $R$ meromorphic on $\mathbb{D}\times\mathbb{P}^1$, satisfying that $R^{-1}(0)$ contains a component that is tangent to the central fiber $E=\{0\}\times\mathbb{P}^1$ at the point $(0,0)$. By the genericity condition on $R$ we have that $R_0(t)=R(0,t)$ has a simple critical point at $0$.  Remark that if we post-compose $R(x,t)$ with a non-constant rational function $Q$ on $\mathbb{P}^1$, the level sets of $Q\circ R$ still define the same foliation. By choosing $Q$ and the $x$ coordinate appropriately we claim that we can suppose that the first integral $R$ for $\mathcal{F}_1$ satisfies

\begin{enumerate}
\item \label{item: critical set} For any critical {\em value} $v\neq 0$ of $R_0$ except possibly for one of them,  there is a connected component of $R^{-1}(v)$ that is not critical for $R$, intersecting $0\times\mathbb{P}^1$ on two points $q,h(q)$ where $h$ is the involution associated to $\mathcal{F}_1$ at $(0,0)$.

  \item \label{item: critical values} $(x,0)\in\mathbb{D}\times\mathbb{P}^1$ is a simple critical point of $R_x(t)=R(x,t)$ with critical value $R_x(0)=x$.

\end{enumerate}
To prove that condition \ref{item: critical set}. can be attained, take a domain $D$ where $h:D\rightarrow D$ is conjugated to a rotation and each leaf cutting $D\setminus 0$ is a disc intersecting $D$ on two points. Take a round disc $D_r\subset R_0(D)$ containing $0$. By composing $R_0$ with a Moebius transformation we can suppose that $D_r=\mathbb{H}$, the upper half plane in $\mathbb{C}$, and the critical values $v_1,\ldots,v_k\in\mathbb{C}\setminus R_0(D)$ of $R_0$  belong to a small neighbourhood of $\infty$.
%\begin{figure}[ht!]
%\centering
%\includegraphics[width=50mm]{graph.pdf}
%\caption{Real graph of the quintic p}
%\label{fig:quintica}
%\end{figure}

Next take a polynomial $Q(z)=z^5+a_4z^4+\ldots+a_1z+a_0$ with real coefficients $a_i\in\mathbb{R}$ satisfying that its four critical points $c_1<c_2<c_3<c_4$ in $\mathbb{C}$ lie in $\mathbb{R}$, and the equation $Q(z)=Q(c_i)$ has precisely two distinct real roots  for each $i=1,\ldots,4$. By construction the other two roots of each such equation are complex conjugate. In particular all finite critical values of $Q$ are attained at regular points in $\mathbb{H}$. To show that the finite critical values  of $Q\circ R_0$ are also attained in $D$ it suffices to remark that in a neighbourhood $U_{\rho}=\{z\in\mathbb{H}:|z|>\rho\}$ for $\rho$ sufficiently big $p$ acts like $z\mapsto z^5$ and thus $Q(U_{\rho})$ covers a pointed neighbourhood of infinity. As $Q(v_i)$ are close to $\infty$ we have that $Q(v_i)\subset Q(\mathbb{H})$.

Once condition \ref{item: critical set}  is satisfied, condition \ref{item: critical values} can be obtained by a change of variables. Indeed, if $R$ already satisfies \ref{item: critical set} then in some connected and simply connected neighbourhood $U\subset 0\times\mathbb{P}^1$ where the involution associated to $\mathcal{F}_1$ is defined, we can define two degree two branched coverings whose fibers coincide: on the one hand $R_{|U}$ and on the other, the projection $\pi$ along the leaves of the foliation from $U$ to the set $t=0$. Hence we can parametrize a neighbourhood of $0$ in $t=0$ by the values $x\in W=R(U)$ and then on $W\times\mathbb{P}^1$ the map $(x,t)\mapsto \widetilde{R}(x,t):=R_{\pi\circ R^{-1}(x)}(t)$ is a holomorphic family of rational functions and satisfies $\widetilde{R}(x,0)=x$. By construction $W$ contains all the critical values of $\widetilde{R}$ except maybe for one. The foliation defined by the levels of $\widetilde{R}$ is equivalent to $\mathcal{F}_1$. After applying the Riemann mapping theorem to $W$ we can suppose $W=\mathbb{D}$.

Let $C=\{v_1,\ldots, v_k\}\subset\mathbb{P}^1\setminus 0$ be the set of critical values of $R_0$ different from $0$. By construction, for each $x\in\mathbb{D}\setminus (C\cap\mathbb{D})$ the rational function $R_x$ has degree $d$ and has critical values at $\{x\}\cup C$. Indeed, since the tangency point is simple and unique, there is a unique component of the tangency divisor between $\mathcal{F}_1$ and the vertical fibration, and it corresponds to the set $t=0$ by construction. Each other critical value of $R_0$ produces a critical value of $R_x$ at the point of intersection of the corresponding leaf with the fibre $\{x\}\times\mathbb{P}^1$.  The restriction of $R_x$ to $R_x^{-1}(\mathbb{P}^1\setminus C\cup \{x\})$ defines a topological degree $d$ covernig having monodromy in a conjugacy class of a subgroup $G_x$ of the symmetric subgroup in $d$ symbols. By continuity the class of $G_x$ is constant $G$ for all $x\in\mathbb{D}\setminus C$.  By connectedness of the covering we know that $G$ acts transitively on each fibre.

Let $\mathcal{H}$ be the Hurwitz space associated to the triple $(d,k+1,G)$, that is, the space of isomorphism classes of topological coverings of the sphere minus $k+1$ points having degree $d$ and monodormy conjugated to $G$. Two coverings $X,X^{\prime}$ are isomorphic if there exists a homeomorphism between the covering spaces $H:X\rightarrow X^{\prime}$ such that $\pi=\pi^{\prime}\circ H$, where $\pi,\pi^{\prime}$ denote the covering projections. In particular for two coverings to be equivalent they need to omit the same set of values on the sphere. Let $\mathcal{V}$ be the set of unordered $(k+1)$-uples of distinct points in $\mathbb{P}^1$.
%Let $\Delta\subset (\mathbb{P}^1)^k$ be the discriminant variety. The projection of $(\mathbb{P}^1)^{k}\setminus \Delta$ by $(\mathbb{P}^1)^k\mapsto (\mathbb{P}^1)^k/S_k$ where $S_k$ denotes the symmetric group on $k$ symbols acting by permuting the coordinates, defines an open subset $V$ in the quotient space.
 Hurwitz (see \cite{H} or \cite{F}) showed that the projection $P:\mathcal{H}\rightarrow \mathcal{V}$, defined by associating to any class of coverings the set of values it omits on the sphere, is itself a topological covering map. We have a natural, continuous, non-constant map $f:\mathbb{D}\setminus C\rightarrow \mathcal{V}$ defined by $f(x)=P([R_x])$. If we take the coordinates in $\mathbb{P}^1$ we took before it can be written as $f(x)=[\{x,v_1,\ldots, v_k\}]\in \mathcal{V}$ and it extends naturally to a map $f:\mathbb{P}^1\setminus C\rightarrow \mathcal{V}$ that is actually holomorphic. To lift $f$ to a map $F: \mathbb{P}^1\setminus C\rightarrow \mathcal{H}$ continuously it suffices to guarantee that at the fundamental group level we have the inclusion $\text{Im}f_*\subset \text{Im} P_*$. This condition is satisfied since we can find generators $\gamma_1,\ldots,\gamma_{k-1}$ of the fundamental group of $\mathbb{P}^1\setminus C$ whose images lie in $\mathbb{D}\setminus (C\cap\mathbb{D})$, and thus the loops $t\mapsto f(\gamma_i(t))$ in $\mathcal{V}$ lift to loops $t\mapsto [R_{\gamma_i(t)}]$ in $\mathcal{H}$. The resulting $F$ has finite fibers and is holomorphic when we consider  the unique complex structure on $\mathcal{H}$ for which $P$ is holomorphic (recall that $\mathcal{V}$ already carries a holomorphic structure).

 For each $x\in \mathbb{P}^1\setminus C$, $F(x)$ defines a unique branched covering of the whole sphere, and hence can be considered as a rational function of degree $d$. The new holomorphic map $\mathbb{P}^1\setminus C\rightarrow \text{Rat}_{\leq d}$ so defined is holomorphic and has finite fibres. Hence it has no essential singularity and it extends to a holomorphic map ${\bf F}: \mathbb{P}^1\rightarrow \text{Rat}_{\leq d}$.

  By construction and uniqueness of complex structure on the sphere, there exists for each $x\in \mathbb{D}\setminus C$ a Moebius transformation $H_x$ such that $R_x\circ H_x=F(x)$. In particular, by pulling $R$ back by  the change of coordinates $(x,t)\mapsto (x,H_x(t))$ defined in a neighbourhood of $x=0$ we have that the the germ of $x\mapsto {\bf F}(x)$ at $0$ describes the pull back of the foliation $\mathcal{F}_1$. This foliation extends to $\mathbb{P}^1\times\mathbb{P}^1$ by the level sets of $F(x,t)={\bf F}(x)(t)$. By blowing up a point of transversality of the foliation and the central fibre and contracting the strict transform of the fibre we obtain  a foliation in $\mathbb{P}^1\times \mathbb{P}^{1}$ having a singularity in $T_1$ with the same holonomy involution as $\mathcal{G}$ modulo conjugation by the Moebius transforamtion $H_0$. By the Remark of Section \ref{adapt} two foliations in $T_1$ having the same involution modulo conjugation by a Moebius transformation are analytically equivalent. Hence we have that the germ $\mathcal{G}$ is equivalent to the germ of that singularity. The obtained foliation is obviously defined by polynomial equations.
\vspace{0.2cm}

This proof cannot be extended to other foliations in $T$ in general because there appear many components of the tangency divisor between ${\cal F}_1$ and the vertical foliation  and there is no way of finding a coordinate where all the curves of critical values can be extended in the same parametrization to $\mathbb{P}^1$. Even if the extension existed there would be intersections of the parametrized curves of critical values and we would have no control over the monodromies around those intersection points. It is for this reason that we will change our point of view and, instead of trying to extend the germ of curve $x\mapsto R_x$ we will try to approximate it by one that extends and use that for good approximations the associated foliations are locally equivalent.

\subsection {Critical Points}

We start by analysing the curves of critical points of $x \mapsto R_x(t)$. We will assume that there exists a fixed neighborhood $U$ (independent of $x$) of $\infty \in {\mathbb P}^1$ such that no critical point is inside this neighborhood. We have the following possibilities:

\begin {itemize}

\item [(i)] the leaf of ${\mathcal F}_1$ that passes through a critical point of $R_0(t)$ (of order $m \in {\mathbb N}$) is transversal to $\{0\}\times {\mathbb P}^1$; we parametrise the leaf as $x\mapsto (x,f(x))$. Since the first integral assumes a constant value along each nearby leaf, we see that each point $(x,f(x))$ is also a critical point of order $m$ of $R_x(t)$. Consequently the curve $t-f(x)$ is contained in the singular set of the foliation defined by $d\,R=0$; we call such a curve of critical points (or singular points) a {\it level type curve}. We may write locally (assuming $t_0=0$ for simplicity) that
$$
R(x,t)= a + (t-f(x))^{m+1}h(x,t)
$$
where $a\in {\mathbb C}$, $h(0,0)\ne 0$. Therefore
$$
d\,R= [(m+1)(t-f(x))^{m}h+ (t-f(x))^{m+1}\dfrac{\partial h}{\partial x}]d\,x
$$
$$
+[-(m+1)(t-f(x))^{m}hg^{\prime} + (t-f(x))^{m+1}\dfrac{\partial h}{\partial t}]d\,t
$$
The foliation $d\,R=0$ has $(t-f(x))^{m}$ as its equation of zeroes. The equation of ${\mathcal F}_1$ is then $\frac {dR}{(t-f(x))^{m}}=0$.

\item [(ii)] the critical point $(0,t_0)$ is a point of tangency of ${\mathcal F}_1$ with $\{0\}\times {\mathbb P}^1$; it gives rise to a curve  of critical points of $R_x(t)$, or points of tangency between ${\mathcal F}_1$ and the vertical lines $x=const$, which crosses $\{0\}\times {\mathbb P}^1$ at the point $(0,t_0)$ (we put again $t_0=0$). The foliation ${\mathcal F}_1$ is obtained in a neighborhood of $(0,0)$ once we divide $d\, R=0$
by the equation of its zeroes. If a component of the curve of critical points is invariant by ${\mathcal F}_1$, it necessarily coincides with the leaf which is tangent to  $\{0\}\times {\mathbb P}^1$ at $(0,0)$; we call it also a {\it level type curve}  of critical points (of some order $M$). It has as equation $x-g(t)=0$, where $g(t)= t^{l+1}{\tilde g}(t)$ with $l\ge 1$ and ${\tilde g}(0)\neq 0$. We apply the same argument as in case (i) to a neighborhood of a point
of this curve for which $x\ne 0$ and conclude that $(x-g(t))^{M}=0$ is inside the set of zeroes of $d\,R$ (a fortiori in a neighborhood of $(0,0)$ as well).

Now let us analyse the case of a component of a {\it non-invariant curve} of critical points, that is, one that is not ${\mathcal F}_1$-invariant. We observe that the zeroes of $d\,R$ are inside the zeroes of $\frac {\partial R}{\partial t}=0$. Locally at a point where $x\ne 0$ we have
$$
R(x,t)= a(x)+ (t-u(x))^{l+1}h(x,t)
$$
where $a(x)$ is not constant (otherwise we would have case (i)), $h(0,0)\ne 0$, $l\ge 1$ and $t-u(x)=0$ is the local equation of the component.  It follows from
$$
d\,R= [{a^{\prime}}(x)-(l+1){(t-u(x))^l}{f^{\prime}}(x)h+{(t-u(x))^{l+1}}\dfrac{\partial h}{\partial x}]d\,x
$$
$$
+[(l+1){(t-u(x))^l}h+ {(t-u(x))^{l+1}}\dfrac{\partial h}{\partial t}]d\,t
$$

that the coefficients of $d\,x$ and $d\,t$ have no common factors; therefore there is no new curve of zeroes arising from the type of curve of critical points under consideration. We conclude that in a neighborhood of $(0,0)$ we just have to take   $\frac {dR}{(x-g(t)))^{M}}=0$ in order to define ${\mathcal F}_1$. Of course it may happen that the leaf of ${\mathcal F}_1$ which is tangent to $\{0\}\times {\mathbb P}^1$ is not a level type curve of critical points.

\end {itemize}

We may summarise this information about the zeroes of $d\,R$ as follows:

\begin {itemize}

\item there are curves of level type $x\mapsto f_1(x),\dots,f_k(x)$ which correspond to critical points of orders $m_1,\dots, m_k$ ; these curves are transversal to $\{0\}\times {\mathbb P}^1$, and locally $R(x,t)=a_j +(t-f_j(x))^{m_j+1}h_j(x,t)$. Locally at each of these critical points the foliation ${\mathcal F}_1$ is given by the equation $\frac {dR}{(t-f_j(x))^{m_j}}=0$.

\item there are $(l_1+1),\dots,(l_s+1)$\,-valued curves $x\mapsto P_1(x),\dots,P_s(x)$ of critical points of orders $M_1 ,\dots,M_s $ which are curves of level type (for $x\ne 0$); each curve $P_j^{\prime}= \cup_{x}P_j(x)$ is tangent    to $\{0\}\times {\mathbb P}^1$ in order $l_j$ at a critical point of $R_0$; its equation is $x-g_j(t)=0$ with  $g_j(t)= t^{l_j+1}{\tilde g}_j(t)$,  $l\ge 1$ and ${\tilde g}_j(0)\neq 0$. Locally at each of these critical points the foliation ${\mathcal F}_1$ is given by the equation $\frac {dR}{(x-g_j(t))^{M_j}}=0$; we have also $R\equiv A_j$ along each curve $x-g_j(t)=0$.

\end {itemize}

\subsection {Proof of Theorem 3}

Now we will consider the algebraic variety  which is the closure of the space of degree $d$ rational functions of ${\mathbb P}^1$ which have the configuration of critical points we presented, namely:

\begin {enumerate}

\item[*] the rational function has values $a_1,\dots,a_k$ at  critical points which have orders $m_1,\dots,m_k$ respectively.

\item[**] the rational function  has values $A_1,\dots,A_s$ at $(l_1+1),\dots, (l_s+1)$ critical points which have orders $M_1,\dots,M_s$ respectively.

\end {enumerate}
%\vspace {0.1in}
Let us denote also by $R$ the curve given by $R(x)=R_x$; it belongs to a smooth stratum $B$ of this variety for $x\ne 0$ small and $R(0)$ belongs to $\bar B$, which is also an algebraic variety. Let $\pi$ be a desingularisation of $\bar B$ and of $R$ at the point $R(0)$. The strict transform ${\tilde R}$ of $R$ crosses the boundary of $\pi^{-1}(B)$ at a smooth point ${\tilde R}(0)\in \pi^{-1}(\bar B)$.
%We remark that the foliation $\mathcal F_1$ of ${\mathbb D}\times {\mathbb P}^{1}$ can also %be considered as a foliation of ${\tilde R}\times {\mathbb P}^{1}$ in a obvious way.
We have a foliation in  ${\tilde R}\times {\mathbb P}^{1}$ given by the level curves of the
meromorphic function $(\tilde p,t) \mapsto R_{\pi (\tilde p)}(t)$, which is conjugated to   $\mathcal F_1$ (since $x\in {\mathbb D} \mapsto R_x$ is injective, it works as a desingularization of the curve $R$).

 Next we take an algebraic curve $\tilde S$ in $\pi^{-1}(\bar B)$ which passes through ${\tilde R}(0)$ smoothly with order of tangency $N$ as big as we wish with $\tilde R$ at the point ${\tilde R}(0)$; the choice of $N$ will depend on the statements which will follow. Consequently in $\bar B$ we may choose an algebraic family $S$ of rational functions  parametrized by a map of $x\in {\mathbb D}$ near the point $S(0)=R(0)$ such that both associated foliations $d\,R=0$ and $d\,S=0$ are as close as we wish in ${\mathbb D}\times {\mathbb P}^1$ (in fact, we need to cover ${\mathbb D}\times {\mathbb P}^1$ by two coordinates systems; in the chart that contains $\{0\}\times \{\infty\}$  we use $R=const$ and $S=const$ to define the associated foliations, which are both regular ones; in the chart that contains $\{0\}\times  \{0\}$ the foliations $d\,R=0$ and $d\,S=0$ are singular).

Next we need to prove that after eliminating the singularities of $d\,S=0$ we obtain a foliation
which is regular and has the same type of tangencies with $\{0\}\times {\mathbb P}^1$ as $\mathcal F_1$.

Let us fix a family of disjoint polydiscs, one for each critical point of $R_0=S_0$. If $(0,t_j)$ is a critical point, we take ${\Delta}_j=\{(x,t); |x|\leq \epsilon, |t-t_j|\leq \epsilon\}$. If $S_x$ is sufficiently close to $R_x$ and $\epsilon$ is small, the configuration of critical points of $S_x$ in each set $K_j=\{(x,t); \frac {\epsilon}{2} \leq |x|\leq \epsilon, |t-t_j|\leq \epsilon\}$ is the same as the configuration of $R_x$. This means that for $S_x$ we have in each $K_j$:

\begin {itemize}

\item there are new curves of level type $x\mapsto {\hat f}_1(x),\dots,{\hat f}_k(x)$ which correspond to critical point of orders $m_1,\dots,m_k$; $S$ takes the values $a_1,\dots,a_k$ along these curves.

\item there are new $(l_1 +1),\dots,(l_s +1)$-valued curves $x\mapsto {\hat P}_1,\dots {\hat P}_s(x)$ which are curves of level type corresponding to critical points of orders $M_1,\dots,M_s$; $S$ takes the values $A_1,\dots,A_s$ along these curves.
\end {itemize}

\noindent Since the set of critical points of $S_x$ inside each ${\Delta}_j$ is an analytic curve, we conclude that the critical curve of level type that lies in $D_j$ has an extension which passes through the point $(0,t_j)$ and reproduces the same type of the corresponding critical curve of $R_x$. The singular set of the foliations $d\,R=0$ and $d\,S=0$ (which appear along the curves of critical points of level type because of condition **) are then sufficiently close (if $R_x$ and $S_x$ are sufficiently close) to allow us to conclude that the regular foliations associated to them (after elimination of singularities) are also close. Notice that the equality $R_0=S_0$ implies that the involutions associted to both foliations coincide.

Now we blow-up the point $(0,\infty)\in {\mathbb D}\times {\mathbb P}^1$; we obtain  two foliations (one for $R_x$ and the other one for $S_x$) defined in  neighborhoods of the strict transform $E^{\prime}_1$ of $\{0\}\times \mathbb P^1$ that  have the same holonomy invariants at the points of tangency with $E^{\prime}_1$ since $R_0=S_0$ . Furthermore, they may be assumed to be close enough as to allow application of \cite {CAL}; consequently they are conjugated.

We remark that we have a global foliation associated to the strict transform   $\tilde S$ of $S$ (the closure of   $\pi^{-1}(S\setminus \{S(0)\})$). We define in ${\tilde S}\times {\mathbb P^1}$ the meromorphic function which
at a fiber $\{p\}\times {\mathbb P^1}$ is exactly ${\tilde S}_p$; we have then inside the algebraic variety ${\tilde S}\times {\mathbb P^1}$ the algebraic foliation given by the level curves of this meromorphic function (it is easy to see that this foliation is equivalent to $dS=0$ in a neighborhood of $\pi^{-1}(S(0))\times \mathbb P^1$. Now we just blow up ${\tilde S}\times {\mathbb P^1}$ at the point $(0,\infty)\in {\mathbb D}\times {\mathbb P}^1$, and conclude by blowing down the exceptional divisor $E_1^{\prime}$.

\bigskip
\begin{minipage}{0.49\linewidth}
G. Calsamiglia\\
{\scriptsize Instituto de Matem\'{a}tica e Estat\'\i stica \\
Universidade Federal Fluminense\\
Rua M\' ario Santos Braga s/n\\
24020-140, Niter\'{o}i, Brazil \\
gabriel@mat.uff.br }
\end{minipage}
\begin{minipage}{0.49\linewidth}
P. Sad\\
{\scriptsize IMPA\\
Estrada Dona Castorina 110\\
Rio de Janeiro / Brasil 22460-320\\
sad@impa.br}
\end{minipage}

\end{document}